\documentclass[12pt]{amsart}

\usepackage{amsthm}
\usepackage{amsmath}
\usepackage{amssymb}
\usepackage{eucal}
\usepackage{amsfonts}
\usepackage{amsmath}
\usepackage[Lenny]{fncychap}
\usepackage{enumerate}
\usepackage{amssymb}
\usepackage[english]{babel}
\usepackage[breaklinks=true]{hyperref}
\usepackage[all]{xy}
\usepackage{fancybox}
\usepackage{latexsym,mathrsfs,longtable,lscape,dsfont,yfonts}
\usepackage{srcltx}
\usepackage{soul}

\ifx\pdfpagewidth\undefined 
\usepackage[T1]{fontenc}
\else 
\usepackage[OT1]{fontenc} 
\fi 
\usepackage{amsfonts,amssymb,latexsym,longtable,amsmath}

\newtheorem{thm}{Theorem}[section]
\newcommand{\bthm}{\begin{thm}} \newcommand{\ethm}{\end{thm}}
\newtheorem{prop}[thm]{Proposition}
\newcommand{\bprp}{\begin{prop}} \newcommand{\eprp}{\end{prop}}
\newtheorem{fact}[thm]{Fact}
\newcommand{\bfct}{\begin{fact}} \newcommand{\efct}{\end{fact}}
\newtheorem{prob}[thm]{Problem}
\newcommand{\bprb}{\begin{prob}} \newcommand{\eprb}{\end{prob}}
\newtheorem{quest}[thm]{Question}
\newcommand{\bqtn}{\begin{quest}} \newcommand{\eqtn}{\end{quest}}
\newtheorem{lem}[thm]{Lemma}
\newcommand{\blem}{\begin{lem}} \newcommand{\elem}{\end{lem}}
\newtheorem{claim}[thm]{Claim}
\newcommand{\bclm}{\begin{claim}} \newcommand{\eclm}{\end{claim}}
\newtheorem{cor}[thm]{Corollary}
\newcommand{\bcor}{\begin{cor}} \newcommand{\ecor}{\end{cor}}
\newtheorem{conj}[thm]{Conjecture}
\newcommand{\bcnj}{\begin{conj}} \newcommand{\ecnj}{\end{conj}}
\theoremstyle{definition}
\newtheorem{defn}[thm]{Definition}
\newcommand{\bdfn}{\begin{defn}} \newcommand{\edfn}{\end{defn}}
\newtheorem{spec}[thm]{Specializing}
\newcommand{\bspc}{\begin{spec}} \newcommand{\espc}{\end{spec}}
\theoremstyle{remark}
\newtheorem{rem}[thm]{Remark}
\newcommand{\brem}{\begin{rem}} \newcommand{\erem}{\end{rem}}
\newtheorem{cnv}[thm]{Convention}
\newcommand{\bcnv}{\begin{cnv}} \newcommand{\ecnv}{\end{cnv}}
\newtheorem{exam}[thm]{Example}
\newcommand{\bexm}{\begin{exam}} \newcommand{\eexm}{\end{exam}}
\newcommand{\bpf}{\begin{proof}} \newcommand{\epf}{\end{proof}}

\newtheorem{thmy}{\textbf{Theorem}}
\newenvironment{thmx}{\stepcounter{thm}\begin{thmy}}{\end{thmy}}

\newcommand{\R}{\mathbb R}
\newcommand{\C}{\mathbb C}
\newcommand{\Z}{\mathbb Z}

\newcommand{\N}{\mathbb N}

\newcommand{\bD}{\mathbb D}

\renewcommand{\phi}{\varphi}
\renewcommand{\theta}{\vartheta}

\newcommand{\w}{{\rm w}}

\newcommand{\U}{\mathbb{U}}

\newcommand{\mkp}{\medskip}

\def\defi{\buildrel\rm def \over=}

\begin{document}

\title[A dichotomy property for locally compact groups]{A dichotomy property for locally compact groups}

\author[M. Ferrer]{Mar\'ia V. Ferrer}
\address{Universitat Jaume I, Instituto de Matem\'aticas de Castell\'on,
Campus de Riu Sec, 12071 Castell\'{o}n, Spain.}
\email{mferrer@mat.uji.es}

\author[S. Hern\'andez]{Salvador Hern\'andez}
\address{Universitat Jaume I, Departamento de Matem\'{a}ticas,
Campus de Riu Sec, 12071 Castell\'{o}n, Spain.}
\email{hernande@mat.uji.es}

\author[L. T\'arrega]{Luis T\'arrega}
\address{Universitat Jaume I, IMAC and Departamento de Matem\'{a}ticas,
Campus de Riu Sec, 12071 Castell\'{o}n, Spain.}
\email{ltarrega@uji.es}

\thanks{ Research Partially supported by the Spanish Ministerio de Econom\'{i}a y Competitividad,
grant MTM2016-77143-P (AEI/FEDER, EU),
and the Universitat Jaume I, grant P1171B2015-77. The second
author also acknowledges partial support by Generalitat Valenciana,
grant code: PROMETEO/2014/062.}

\begin{abstract}
We extend to metrizable locally compact groups Rosenthal's theorem describing those Banach spaces containing no copy of $\ell_1$.
For that purpose, we transfer to general locally compact groups the notion of interpolation (\emph{$I_0$}) set, which was defined
by Hartman and Ryll-Nardzewsky \cite{Hartman1964} for locally compact abelian groups.
Thus we prove that for every sequence $\{g_n\}_{n < \omega}$ in a locally compact group $G$, then either $\{g_n\}_{n < \omega}$ has a weak Cauchy subsequence
or contains a subsequence that is an $I_0$ set. This result is subsequently applied to obtain sufficient conditions
for the existence of Sidon sets in a locally compact group $G$, an old question that remains open since 1974
(see \cite{Lopez1975} and \cite{Figa-Talamanca1977}).
Finally, we show that every locally compact group strongly respects compactness extending thereby a result by Comfort, Trigos-Arrieta, and Wu
\cite{Comfort1993}, who established this property for abelian locally compact groups.

\end{abstract}

\thanks{{\em 2010 Mathematics Subject Classification.} Primary 22D05; 43A46. Secondary 22D10; 22D35; 43A40; 54H11\\
{\em Key Words and Phrases:} Locally compact group; $I_0$ set; Sidon set; Interpolation set; Weak topology; Bohr compactification; Eberlein compactification}

\dedicatory{Respectfully dedicated to Professor Wis Comfort}

\date{\today}

\maketitle \setlength{\baselineskip}{24pt}

\section{Introduction}

A well-known result of Rosenthal establishes that if $\{x_n\}_{n<\omega}$ is
a bounded sequence in a Banach space $X$, then either $\{x_n\}_{n<\omega}$ has a weak
Cauchy subsequence, or $\{x_n\}_{n<\omega}$ has a subsequence equivalent to the
usual $\ell_1$-basis (and then $X$ contains a copy of $\ell_1$).
In this paper,we look at this result for locally compact groups.
More precisely, our main goal is to extend Rosenthal's dichotomy theorem on Banach spaces to
locally compact groups and their weak topologies. First, we need some definitions and
basic results.

Given a locally compact group $(G,\tau)$, we denote by $Irr(G)$ the set of all continuous
unitary irreducible representations $\sigma$ defined on $G$.
That is, continuous in the sense that each matrix coefficient function
$g\mapsto \langle\sigma(g)u,v\rangle$ is a continuous map of $G$ into the complex plane.
Thus, fixed $\sigma\in Irr(G)$, if $\mathcal{H}^{\sigma}$ denotes the Hilbert space associated to $\sigma$,
we equip the unitary group $\U(\mathcal H^\sigma)$ with the weak (equivalently, strong) operator topology.
For two elements $\pi$ and $\sigma$ of $Irr(G)$, we write $\pi\sim\sigma$
to denote the relation of unitary equivalence and we denote by $\widehat G$ the \emph{dual object} of $G$,
which is defined as the set of equivalence classes in ($Irr(G)/{\sim}$).
We refer to \cite{Dixmier1964,Bekka2008} for all undefined
notions concerning the unitary representations of locally compact groups.

Adopting, the terminology introduced by Ernest in \cite{Ernest1971}, set $\mathcal{H}_n\defi \C^n$ for $n=1,2, \ldots$;
and $\mathcal H_{0}\defi l^2(\Z)$. The symbol $Irr^C_n(G)$ will denote the set of irreducible unitary representations of $G$
on $\mathcal{H}_n$, where it is assumed that every set $Irr^C_n(G)$ is equipped with the compact open topology.
Finally, define $Irr^C(G)=\bigsqcup\limits_{n\geq 0} Irr^C_n(G)$ (the disjoint topological sum).

We denote by $G^w=(G,w(G,Irr(G))$ (resp. $G^{w_C}=(G,w(G,Irr^C(G)))$) the group $G$ equipped with
the weak (group) topology generated by $Irr(G)$ (resp. $Irr^C(G)$).
Since equivalent representations define the same topology, we have $G^w=(G,w(G,\widehat G))$. That is, the \emph{weak topology}
is the initial topology on $G$ defined by the dual object. Moreover, in case $G$ is a separable, metric, locally compact group,
then every irreducible unitary representation acts on a separable Hilbert space and, as a consequence, is unitary equivalent
to a member of $Irr^C(G)$. Thus $G^w=(G,w(G,Irr^C(G)))=G^{w_C}$ for separable, metric, locally compact groups.
We will make use of this fact in order to avoid
the proliferation of isometries (see \cite{Dixmier1964}).
In case the group $G$ is abelian,
the dual object $\widehat G$ is a group, which is called \emph{dual group},
and the weak topology of $G$ reduces to the weak topology generated by
all continuous homomorphisms of $G$ into the unit circle $\mathbb T$. That is,
the weak topology coincides with the
so-called \emph{Bohr topology} of $G$, that we recall next for the reader's sake.

With every (not necessarily abelian) topological group $ G $ there is associated a compact Hausdorff group
$ bG $, the so-called \emph{Bohr compactification} of $ G $, and a continuous homomorphism
$ b $ of $ G $ onto a dense subgroup of $ bG $ such that $ bG $ is characterized by the following universal property:
given any continuous homomorphism $ h $ of $ G $ into a compact group $ K $,
there is always a continuous homomorphism $ \bar {h} $ of $ bG $ into $ K $ such that
$ h = \bar {h} \circ b $ (see \cite [V \S 4]{Heyer1970}, where a detailed study on $bG$ and their properties is given).
In Anzai and Kakutani \cite{Anzai1943} $bG$ is built when
$ G $ is locally compact abelian (LCA). However, most authors agree that it was A. Weil \cite{Weil1965} the first
to build $ bG $. Weil called $ bG $
``Groupe compact attach\'{e} \`{a} $ G $ ''. The name of Bohr compactification was given by
  Alfsen and Holm \cite{Alfsen1962} in the context of arbitrary topological  groups.
The \emph{Bohr topology} of a topological group $G$ is the one that inherits as a subgroup of $bG$.

The weak topology of a topological group plays a role analogous to the weak topology in a Banach space.
Therefore, it is often studied in connection to the original topology of the group.
For instance, it can be said that the preservation
of compact-like properties from $G^{w}$ to $G$ concerns
``uniform boundedness'' results and, in many cases, it can
be applied to prove the continuity of
 certain related algebraic homomorphisms.

Our main result establishes that for every sequence $\{g_n\}_{n < \omega}$ in a locally compact group $G$,
then either $\{g_n\}_{n < \omega}$ has a weak Cauchy subsequence
or contains a subsequence that is an $I_0$ set.
As a consequence, we obtain some sufficient conditions for the existence of (weak) Sidon sets in locally compact groups.
It is still an open question whether every infinite subset of a locally compact group $G$ contains a (weak) Sidon subset
(see \cite{Lopez1975,Figa-Talamanca1977}).

As far as we know, the first result about the existence of $I_0$ sets was given by Hartman and Ryll-Nardzewski \cite{Hartman1964},
who considered the weak topology associated to a locally compact abelian (LCA, for short) group and introduced the notion
of \emph{interpolation} (or $I_0$) set.
As they defined it, a subset $E$ of a LCA group $G$
is an $I_0$ set if  every bounded function on E is the restriction of an almost \emph{periodic function}  on $G$
(here, it is said that a complex-valued
function $f$ defined on $G$ is \emph{almost periodic} when is the restriction of a continuous function defined on $bG$).
Alternatively, one can define this notion without recurring to the Bohr compactification using the Fourier transform.
Thus $E$ is an $I_0$ set if every bounded function on E is the restriction of the Fourier transform of a discrete measure on $G$.

Therefore, an $I_0$ set is a subset $E$ of $G$ such that any bounded map on $E$ can be interpolated by a continuous function on $bG$.
As a consecuence, if $E$ is a countably infinite $I_0$ set, then $\overline{E}^{\, bG}$ is canonically homeomorphic to  $\beta\omega$.
The main result given by Hartman and Ryll-Nardzewski is the following:

\bthm[\cite{Hartman1964}] \label{douw}
Every LCA group $G$ contains $I_0$ sets.
\ethm
\mkp

For the particular case of discrete abelian groups, van Douwen achieved a remarkable progress
by proving the existence of $I_0$ sets in very general situations.
His main result can be formulated in the following way:

\bthm[\cite{Douwen1990}, Theorem 1.1.3] \label{douw}
Let $G$ be a discrete Abelian group and let $A$ be
an infinite subset of $G$. Then there is a subset $B$ of
$A$  with $|B|=|A|$ such that $B$ is
an $I_0$ set.
\ethm
\mkp

In fact, van Douwen
extended his result to the real line but left unresolved
the question for LCA groups. In general,
the weak topology of locally compact groups has been considered by many workers so far, specially for abelian groups,
where the amount of important results in this direction is vast
(see \cite{Graham2013} for a recent a comprehensive source
on the subject).

For those locally compact groups that can be injected in their
Bohr compactification, the so-called \emph{maximally almost pe\-rio\-dic groups},
the existence of $I_0$ sets was clarified in \cite{Galindo_Hernandez2004}.
However, many (non-abelian) locally compact groups cannot be injected in their Bohr compactification
(that can become trivial in some cases).
Thus, the question of extending Rosenthal's result to general locally compact groups remained open until now.
\mkp

The starting point of our research lies on three celebrated results on $C(X)$, the space of continuous functions on
a Polish space $X$, equipped with the pointwise convergence topology.
The first two basic results are Rosenthal's dichotomy theorem \cite{Rosenthal1974a} and a theorem by Bourgain, Fremlin and Talagrand
about compact subsets of Baire class $1$ functions \cite{Bourgain1978a}, that we present in the way
they are formulated by Todor\v{c}evi\'{c} in \cite{Todorcevic1997}.

\bthm\emph{(H. P. Rosenthal)}\label{ros}
If $X$ is a Polish space and $\{f_n\}\subseteq C(X)$ is a pointwise bounded sequence, then either $\{f_n\}$ contains a convergent
subsequence or a subsequence whose closure in $\mathbb R^X$ is homeomorphic to $\beta\omega$.
\ethm
\mkp

\bthm\emph{(J. Bourgain, D.H. Fremlin, M. Talagrand)}\label{bft}
Let $X$ be a Polish space and let $\lbrace f_n\rbrace_{n<\omega}\subseteq C(X)$ be a pointwise bounded sequence.
The following assertions are equivalent (where the closure is taken in $\mathbb R^X$):
\begin{enumerate}[(a)]
\item $\lbrace f_n\rbrace_{n<\omega}$ is sequentially dense in its closure.
\item The closure of $\lbrace f_n\rbrace_{n<\omega}$ contains no copy of $\beta\omega$.
\end{enumerate}
\ethm
\mkp

Our third starting fact is extracted from a result by Pol \cite[p. 34]{Pol1984}, that again was formulated in different terms (cf. \cite{Cascales2000}).
Here, we only use one of the implications established by Pol.

\bthm\emph{(R. Pol)}\label{pol}
Let $X$ be a complete metric space, $G$ a subset of $C(X)$ which is uniformly bounded.
If $\overline G^{\mathbb R^X}\nsubseteq B_1(X)$, then $G$ contains a sequence whose closure in $\mathbb R^X$ is homeomorphic to $\beta\omega$.
\ethm
\mkp

This paper is divided as follows: The first section is introductory in nature, recounting briefly the history of the topic
and some basic results that are needed along the paper. In the second section, we study sets of continuous functions whose pointwise closure is compact and
contained in the space of all Baire class $1$ functions. We analyze the special case where the functions are defined from a polish space $X$ to a metric space $M$.
Rosenthal \cite{Rosenthal1974a}, Bourgain \cite{Bourgain1977}, and
Bourgain, Fremlin and Talagrand \cite{Bourgain1978a} and, in a different direction, Todor\v{c}evi\'{c} \cite{Todorcevic1999}
have extensively studied the compact subsets of $B_1(X)$.
Our aim is to extend some of their fundamental results to the special case where the functions are metric-valued.

The result that connects with the rest of the paper is the extension of the Rosenthal's theorem \cite{Rosenthal1974a}.
Thus, the main goal is a dichotomy-type result for me\-tric locally compact topological groups that has
several applications in the study of general locally compact groups.
For this purpose, we use the notion of $I_0$ set that plays a \emph{r\^ole} anologous to the $\ell_1$-basis in the realm
of locally compact groups.

In the fourth section, we look at Sidon subsets of non-abelian groups and we apply our previous results in order
to obtain sufficient conditions for the existence of these interpolation sets. Here,
a subset $E$ of $G$ is called \emph{(weak) Sidon set} when every bounded function
can be interpolated by a continuous function defined on the Eberlein compactification $eG$ (defined below).
This is a weaker property than the classical notion of \textit{Sidon set}
in general (see \cite{Picardello1973}) but both notions coincide for amenable groups.

To finish the paper we deal with the property of strongly respecting compactness introduced in \cite{Comfort1993}.
Comfort, Trigos-Arrieta and Wu proved that every abelian locally compact group strongly respects compactness.
Having in mind that Hughes \cite{Hughes1972} proved that every locally compact group, not necessarily abelian, respects compactness,
we extend this result verifying that in fact every locally compact group strongly respects compactness.
This improve previous results by Comfort, Trigos-Arrieta and Wu \cite{Comfort1993}  and Galindo and Hern\'{a}ndez \cite{Galindo_Hernandez2004} mentioned above.

We now formulate our main results. In the sequel $\w G$ will denote the weak compactification of a locally compact group $G$
(defined below) and $inv(\w G) \defi \{x\in\w G : xy = yx = 1\ \text{for some}\ y\in\w G\}$ will designate the group of units
of $\w G$.

\begin{thmx}\label{Teo_A}
Let $(G,\tau)$ be a metric locally compact group and let $\lbrace g_n\rbrace_{n<\omega}$ be a sequence in $G$. Then,
either $\lbrace g_n\rbrace_{n<\omega}$ contains a weak Cauchy subsequence or an $I_0$ set.
\end{thmx}

\begin{thmx}\label{Cor_C}
Every non-precompact subset of a locally compact group $G$  whose weak closure is placed in $inv(\w G)$ contains an infinite $I_0$ set.
Furthermore, this $I_0$ set is precompact in the weak completion of $G$.
\end{thmx}


\begin{thmx}\label{Teo_Sidon}
Every non-precompact subset of a locally compact group  whose weak closure is placed in $inv(\w G)$ contains an infinite weak Sidon set.
Again, this Sidon set is precompact in the weak completion of $G$.
\end{thmx}

\begin{thmx}\label{Teo_SRC}
Every locally compact group $G$ strongly respects compactness.
\end{thmx}


\section{Baire class $1$ functions}

Let $X$ be a topological space and let $(M,d)$ be a metric space. We let $M^X$ (resp. $C(X,M)$) denote
the set of functions (resp. continuous functions) from $X$ to $M$,
equipped with the product (equivalently, pointwise convergence) topology, unless other\-wise stated.
We assume without loss of generality that $d(x,y)\leq 1$ for all $x,y$ in $M$.

A function $f:X\rightarrow M$ is said to be \textit{Baire class $1$} if there is a sequence of continuous functions that converges pointwise to $f$. We denote by $B_1(X,M)$ the set of all $M$-valued Baire $1$ functions on $X$. If $M=\R$ we simply write $B_1(X)$.
A compact space $K$ is called \textit{Rosenthal compactum} if $K$ can be embedded in $B_1(X)$
for some Polish space $X$. Let $X$ be a topological space, the \emph{tightness} of $X$, denoted $tg(X)$,
is the smallest infinite cardinal $\kappa$ such that for any subset $A\subseteq X$ and any point $x\in \overline A$ there is a subset $B\subseteq A$
with $\vert B\vert\leq \kappa$ and $x\in\overline{B}$. We write $[A]^{\leq \omega}$ to denote the set of all countable subsets of $A$.\\

Given a subset $G\subseteq C(X,M)$, it defines an equivalence relation on $X$  by $x\sim y$ if and only if $g(x)=g(y)$ for all $g\in G$.
 If $\widetilde{X}=X/{\sim}$ is the quotient space and $p:X\rightarrow \widetilde{X}$ denotes the canonical quotient map,
 each $g\in G$ has associated a map
$\tilde{g}\in C(\widetilde{X},M)$ defined as $\tilde{g}(\tilde{x})\defi g(x)$ for any $x\in X$ with $p(x)=\tilde{x}$.
Furthermore, if $\tilde{G}\defi \{\widetilde{g} : g\in G\}$, we can extend this definition to the pointwise closure of $\tilde{G}$.
Thus, each $g\in\overline{G}^{M^X}$ has associated a map $\tilde{g}\in\overline{\tilde{G}}^{M^{\widetilde{X}}}$ such that $\tilde{g}\circ p=g$. We denote by $X_G$ the topological space $(\widetilde{X},t_p(\tilde{G}))$. Note that $X_G$ is metrizable if $G$ is countable and it is Polish if $X$ is compact and $G$ is countable.

With the terminology introduced above, the following fact is easily verified (see \cite{FerHerTar:T&A}).

\bfct\label{lem_hom}
If $G$ be a subset of $C(X,M)$ such that $\overline{G}^{M^X}$ is compact, then the map
$p^*:M^{X_G} \rightarrow M^X$, defined by $p^*(\tilde{f})=\tilde{f}\circ p$, is
is a homeomorphism of $\overline{\tilde{G}}^{M^{X_G}}$ onto $\overline{G}^{M^{X}}$.
\efct

\subsection{Real-valued Baire class $1$ functions}

A direct consequence of the strong results presented at the Introduction are the following corollaries whose simple proof is included for the reader's sake.

\bcor\label{cor_polish}
If $X$ is a Polish space and $G$ is a uniformly bounded subset of $C(X)$. Then $tg(\overline{G}^{\mathbb R^X})\leq \omega$
if and only if $\overline{G}^{\mathbb R^X}\subseteq B_1(X)$.
\ecor
\bpf
One implication is consequence of a well known result by Bourgain, Fremlin and Talagrand (cf. \cite{Todorcevic1997}). Therefore,
assume that $tg(\overline{G}^{\R^X})\leq~\omega$. If $\overline{G}^{\mathbb R^X}\nsubseteq B_1(X)$, by Theorem \ref{pol}, we can find a sequence $\lbrace {f}_n\rbrace_{n<\omega}\subseteq {G}$ whose closure
in $\mathbb R^{X}$ is canonically homeomorphic to $\beta\omega$.
This implies that $\overline{G}^{\mathbb R^X}$ contains a copy of $\beta \omega$, which is a contradiction.
\epf
\mkp

\bcor\label{cor_polish_2}
Let $X$ be a Polish space, $G$ a uniformly bounded subset of $C(X)$.
The following assertions are equivalent:
\begin{enumerate}[(a)]
\item $tg(\overline{G}^{\mathbb R^X})\leq \omega$.
\item $\overline{{G}}^{\mathbb R^{X}}\subseteq B_1(X)$.
\item $G$ is sequentially dense in $\overline{G}^{\mathbb R^X}$.
\item $\vert \overline{G}^{\mathbb R^X}\vert\leq \mathfrak{c}$.
\item $G$ does not contain any sequence whose closure in $\mathbb R^X$ is homeomorphic to $\beta\omega$.
\end{enumerate}
\ecor
\bpf
$(a)\Leftrightarrow (b)$ is Corollary \ref{cor_polish}.

$(b)\Rightarrow (c)$ was proved by Bourgain, Fremlin and Talagrand (see \cite{Todorcevic1997}).

$(b)\Rightarrow (d)$, $(a)\Rightarrow (e)$ and $(c)\Rightarrow (a)$ are obvious.

$(d)\Rightarrow (b)$ and $(e)\Rightarrow (b)$ are a consequence of Theorem \ref{pol}.

\epf
\mkp


Furthermore, according to results of Rosenthal \cite{Rosenthal1974a} and Talagrand \cite{Talagrand1984} we can also add
the property of  containing a sequence equivalent to the unit basis $\ell_1$.

\bdfn
Let $\lbrace g_n\rbrace_{n<\omega}$ be a uniformly bounded real (or complex) sequence of continuous functions on a set $X$. We say that $\lbrace g_n\rbrace_{n<\omega}$ is equivalent to the unit basis $\ell_1$ if there exists a real constant $C>0$ such that
$$   \sum\limits_{i=1}^N \vert a_i\vert\leq C \cdot\Vert \sum\limits_{i=1}^N a_i g_i\Vert_{\infty}$$
for all scalars $a_1,\ldots,a_N$ and $N\in\omega$.
\edfn

\bthm[Talagrand \cite{Talagrand1984}] \label{talagrand}
Let $X$ be a compact and metric space and $G$ a uniformly bounded subset of $C(X)$.
The following assertions are equivalent:
\begin{enumerate}[(a)]
\item $\overline{{G}}^{\mathbb R^{X}}\subseteq B_1(X)$.
\item Every sequence in $G$ has a weak-Cauchy subsequence.
\item $G$ does not contain any sequence equivalent to the $\ell_1$ basis.
\end{enumerate}
\ethm

\brem
It is pertinent to notice here that using Rosenthal-Dor Theorem \cite{Dor1975},
Talagrand's result formulated above also holds for complex valued continuous functions.
\erem

A slight variation of Corollary \ref{cor_polish}  is also fulfilled if $X$ is a compact space and $G$ is countable.

\bcor\label{cor_1}
Let $X$ be a compact space and $G$ a countable uniformly bounded subset of $C(X)$. Then $tg(\overline{G}^{\mathbb R^X})\leq \omega$
if and only if $G$ does not contain any sequence whose closure in $\mathbb R^X$ is homeomorphic to $\beta\omega$.
\ecor
\bpf
Let $X_G$ be the quotient space associated to $G$ equipped with the topology of pointwise convergence on $G$.
According to Fact \ref{lem_hom}, we may assume  without loss of generality that $X=X_G$ and therefore that is a Polish space.
It now suffices to apply Corollary \ref{cor_polish_2}.
\epf
\mkp

\bcor\label{resultado_1}
Let $X$ be a compact space, $G$ a countable uniformly bounded subset of $C(X)$.
The following assertions are equivalent:
\begin{enumerate}[(a)]
\item $tg(\overline{G}^{\mathbb R^X})\leq \omega$.
\item $G$ does not contain any sequence whose closure in $\mathbb R^X$ is homeomorphic to $\beta\omega$.
\item $\overline{G}^{\mathbb R^X}$ is a Rosenthal compactum.
\item $\vert \overline{G}^{\mathbb R^X}\vert\leq \mathfrak{c}$.
\item $G$ does not contain any subsequence equivalent to the $\ell_1$ basis.
\end{enumerate}
\ecor
\bpf

If $X_G$ denotes the quotient space associated to $G$, then Fact \ref{lem_hom} implies that
$\overline{{G}}^{\mathbb R^{X}}$ is canonically homeomorphic to $\overline{\tilde{G}}^{\mathbb R^{X_G}}$

$(a)\Leftrightarrow (b)$ is Corollary \ref{cor_1}.

$(a)\Rightarrow (c)$ By Corollary \ref{cor_polish}, we have that $\overline{\tilde{G}}^{\mathbb R^{X_G}}\subseteq B_1(X_G)$.
Thus $\overline{\tilde{G}}^{\mathbb R^{X_G}}$ and, consequently, also $\overline{{G}}^{\mathbb R^{X}}$ are Rosenthal compacta.


$(c)\Rightarrow (d)$ and $(d)\Rightarrow (b)$ are obvious.

$(a)\Leftrightarrow (e)$ It follows from Fact \ref{lem_hom}, Corollary \ref{cor_polish} and Theorem \ref{talagrand}.
\epf
\mkp

\bcor\label{cor_1_cor}
Let $X$ be a compact space and $G$ a uniformly bounded subset of $C(X)$. If $tg(\overline{G}^{\mathbb R^X})\leq \omega$,
then $\overline{G}^{\mathbb R^X}\subseteq B_1(X)$.
\ecor
\bpf
Suppose that there is $f\in \overline{G}^{\mathbb R^X}\setminus B_1(X)$. Since $tg(\overline{G}^{\mathbb R^X})\leq \omega$, there is $L\in [G]^{\leq \omega}$ such that $f\in \overline{L}^{\mathbb R^X}$. Therefore $f\in \overline{L}^{\mathbb R^X}\setminus B_1(X)$ and,
by Fact \ref{lem_hom}, we deduce that $\tilde f\in\overline{\tilde{L}}^{\mathbb R^{X_L}}\setminus B_1(X_L)$.
It now suffices to apply Corollaries \ref{cor_polish_2} and \ref{cor_1}.
\epf

The following example shows that Corollary \ref{resultado_1} may fail if one takes $G$ of uncountable cardinality.

\bexm Let $X=[0,\omega_1]$ and set
$$G=\lbrace \chi_{[\alpha,\omega_1]}: \alpha<\omega_1,\ \text{ and}\ \alpha\
\text{it is not a limit ordinal}\rbrace.$$ Then $\vert \overline G^{\mathbb R^X}\vert = \mathfrak{c}$. However
$\chi_{\lbrace\omega_1\rbrace}$ is in the closure of $G$ but 
it does not belong to the closure of any countable subset of $G$.
\eexm

\subsection{Metric-valued Baire class $1$ functions}

The goal in this section is to extend the results obtained for real-valued Baire class $1$ functions
to functions that take value in a metric space. This is accomplished using an idea of Christensen \cite{Christensen1981}.
First, we need the following definition.


Let $(M,d)$ be a metric space that we always assume equipped with a bounded metric. We set
$$ \mathcal{K}\defi\lbrace \alpha:M\rightarrow [-1,1]:\vert \alpha(m_1)-\alpha(m_2)\vert\leq d(m_1,m_2),\quad\forall m_1,m_2\in M\rbrace.$$

\noindent Being pointwise closed and equicontinuous by definition, it follows that $\mathcal{K}$ is a compact and metrizable subspace of $\mathbb R^M$.
For $m_0\in M$, define $\alpha_{m_0}\in \mathbb R^M$ by $\alpha_{m_0}(m)\defi d(m,{m_0})$ for all $m\in M$.
It is easy to check that $\alpha_{m_0}\in \mathcal{K}$.
Given $f\in M^X$ we associate a map
$\check{f}\in \mathbb R^{X\times \mathcal{K}}$ defined by $$\check{f}(x,\alpha)=\alpha(f(x))\ \hbox{for all}\ (x,\alpha)\in X\times \mathcal{K}.$$
In like manner, given any subset $G$ of $M^X$ we set $\check{G}\defi\lbrace \check{f}:f\in G\rbrace$.


\blem\label{lema_K}
Let $X$ be a topological space, $(M,d)$ a metric space and $G\subseteq C(X,M)$. Then:
\begin{enumerate}[(a)]
\item $f\in C(X,M)$ if and only if $\check{f}\in C(X\times \mathcal{K})$.
\item A net $\{g_{\delta}\}_{\delta\in \w }\subseteq C(X,M)$ converges to $f\in M^X$ if and only if
the net $\{\check g_{\delta}\}_{\delta\in \w }\subseteq C(X\times \mathcal K)$ converges to $\check{f}\in \mathbb R^{X\times \mathcal{K}}$.
\item If $\overline{G}^{M^X}$ is compact, then $\overline{G}^{M^X}$ and $\overline{\check{G}}^{\mathbb R^{X\times \mathcal{K}}}$
are canonically homeomorphic.
\end{enumerate}
\elem
\bpf
$(a)$ Suppose that $f\in C(X,M)$ and let $\lbrace (x_{\delta},\alpha_{\delta})\rbrace_{\delta\in \w }\subseteq X\times \mathcal{K}$ be a net that converges to $(x,\alpha)\in X\times \mathcal{K}$.
For every $\delta\in \w $, we have
\begin{eqnarray*}
\vert \alpha_{\delta}(f(x_{\delta}))-\alpha(f(x))\vert &\leq & \vert \alpha_{\delta}(f(x_{\delta}))-\alpha_{\delta}(f(x))\vert+\vert\alpha_{\delta}(f(x))-\alpha(f(x)) \vert\leq\\
    &\leq & d(f(x_{\delta}),f(x))+\vert\alpha_{\delta}(f(x))-\alpha(f(x)) \vert
\end{eqnarray*}

Since $\lbrace f(x_{\delta})\rbrace_{\delta\in \w }$ converges to $f(x)$ and $\lbrace \alpha_{\delta}\rbrace_{\delta\in \w }$ converges to $\alpha$
it follows that
$$\lim\limits_{\delta\in \w }\check{f} (x_{\delta},\alpha_{\delta})= \lim\limits_{\delta\in \w } \alpha_{\delta} (f(x_{\delta}))=\alpha(f(x))=\check{f}(x,\alpha).$$
Conversely, suppose that $\check f\in C(X\times \mathcal K)$ and let $\lbrace x_{\delta}\rbrace_{\delta\in \w }\subseteq X$ a net that converges to $x\in X$. Consider the map $\alpha_{f(x)}\in \mathcal{K}$. We have
$$ \lim\limits_{\delta\in \w } d(f(x_{\delta}),f(x))=\lim\limits_{\delta\in \w } \alpha_{f(x)}(f(x_{\delta}))=\lim\limits_{\delta\in \w } \check{f}(x_{\delta},\alpha_{f(x)})=\check{f}(x,\alpha_{f(x)})=0,$$
That is, $f$ is continuous.

$(b)$ Suppose that $\lbrace g_{\delta}\rbrace_{\delta\in \w }\subseteq C(X,M)$ converges pointwise to $f\in M^X$ and take the associated sequence $\lbrace \check{g}_{\delta}\rbrace_{\delta\in \w }\subseteq C(X\times \mathcal{K})$. Then $\lim\limits_{\delta\in \w }\check{g}_{\delta}(x,\alpha)=\lim\limits_{\delta\in \w }\alpha(g_{\delta}(x))=\alpha(\lim\limits_{\delta\in \w } g_{\delta}(x))=\alpha(f(x))=\check{f}(x,\alpha)$ for all $(x,\alpha)\in X\times \mathcal{K}$.

Conversely, suppose that $\lbrace \check g_{\delta}\rbrace_{\delta\in \w }\subseteq C(X\times \mathcal K)$ converges pointwise to
$\check{f}\in \mathbb R^{X\times \mathcal{K}}$ and let us see that the sequence $\lbrace g_{\delta}\rbrace_{\delta\in \w }\subseteq C(X,M)$ converges
pointwise to $f$.
Indeed, it suffices to notice that for every $x\in X$ and its associated map $\alpha_{f(x)}\in \mathcal K$, we have
$$\vert \check{g}_{\delta}(x,\alpha_{f(x)})-\check{f}(x,\alpha_{f(x)})\vert=\vert \alpha_{f(x)}(g_{\delta}(x))-\alpha_{f(x)}(f(x))\vert=d(g_{\delta}(x),f(x)).$$

$(c)$ Consider the map $\phi:\overline{G}^{M^X}\rightarrow \overline{\check{G}}^{\mathbb R^{X\times \mathcal{K}}}$
defined by $\phi(g)\defi \check{g}$ for all $g\in  \overline G^{\, M^X}$.
By compactness, it is enough to prove that $\phi$ is injective and continuous. The argument verifying the continuity of $\phi$ has been used in (b).
Thus we only verify that $\phi$ is injective.
Assume that $\phi(f)=\phi(g)$ with $f,g\in \overline{G}^{M^X}$, which means $\alpha(f(x))=\alpha(g(x))$ for all $(x,\alpha)\in X\times \mathcal{K}$.
Given $x\in X$, we have the map $\alpha_{g(x)}\in \mathcal{K}$ and, consequently, $\alpha_{g(x)}(f(x))=\alpha_{g(x)}(g(x))=0$ for all $x\in X$. This yields $d(f(x),g(x))=0$ for all $x\in X$, which implies $f=g$. 
\epf

Using the previous lemma we can generalize Theorem \ref{ros} to any metric space. This result will be very useful in the sequel.


\bcor\label{cor_Rosenthal_metrico}
Let $X$ be a Polish space, $(M,d)$ a metric space and $\lbrace f_n\rbrace_{n<\omega}\subseteq C(X,M)$ such that $\overline{\lbrace f_n\rbrace_{n<\omega}}^{M^X}$ is compact.
Then, either $\lbrace f_n\rbrace_{n<\omega}$ contains a pointwise Cauchy subsequence or a subsequence whose closure in $M^X$ is homeomorphic to $\beta\omega$.
\ecor

Now, we are in position of extending, to the setting of metric-valued functions, the results obtained in the previous section for real-valued functions.

\bprp\label{Bourgain_Metrico}
Let $X$ be a Polish space, $(M,d)$ be a metric space and  $G\subseteq C(X,M)$ such that $\overline{G}^{M^X}$ is compact. The following assertions are equivalent:
\begin{enumerate}[(a)]
\item $tg(\overline{G}^{M^X})\leq\omega$
\item $\overline{G}^{M^X}\subseteq B_1(X,M)$.
\item $G$ is sequentially dense in $\overline{G}^{M^X}$.
\item $\vert\overline{G}^{M^X}\vert\leq \mathfrak{c}$.
\end{enumerate}
\eprp
\bpf
If follows from Lemma \ref{lema_K} and Corollary \ref{cor_polish_2}.

%
%
\epf

\bprp\label{BFT_Metrico}
Let $X$ be a Polish space, $(M,d)$ be a metric space and a sequence $\lbrace f_n\rbrace_{n<\omega}\subseteq C(X,M)$ such that $\overline{\lbrace f_n\rbrace_{n<\omega}}^{M^X}$ is compact. The following assertions are equivalent:
\begin{enumerate}[(a)]
\item $\lbrace f_n\rbrace_{n<\omega}$ is sequentially dense in its closure.
\item The closure of $\lbrace f_n\rbrace_{n<\omega}$ contains no copy of $\beta\omega$.
\end{enumerate}
\eprp
\bpf
$(a)\Rightarrow (b)$ is obvious.

$(b)\Rightarrow (a)$ By Lemma \ref{lema_K} we know that  $\overline{\lbrace \check{f}_n\rbrace_{n<\omega}}^{\R^{X\times \mathcal K}}$ contains no copy of $\beta\omega$. Thus by
Theorem \ref{bft}, it follows that $\lbrace \check{f}_n\rbrace_{n<\omega}$ is sequentially dense in its closure. By Lemma \ref{lema_K} we conclude that $\lbrace f_n\rbrace_{n<\omega}$ is sequentially dense in its closure.
\epf

\bcor\label{cor_1m}
Let $X$ be a compact space, $(M,d)$ be a metric space $G\subseteq C(X,M)$ such that $\overline{G}^{M^X}$ is compact. If $tg(\overline{G}^{M^X})\leq \omega$, then $\overline{G}^{M^X}\subseteq B_1(X,M)$.
\ecor
\bpf
It suffices to apply Corollary \ref{cor_1_cor} and Lemma \ref{lema_K}.
\epf
\bcor\label{resultado_2}
Let $X$ be a compact space, $(M,d)$ be a metric space and $G$ a countable subset of $C(X,M)$ such that $\overline{G}^{M^X}$ is compact.
The following assertions are equivalent:
\begin{enumerate}[(a)]
\item $tg(\overline{G}^{M^X})\leq \omega$.
\item $G$ does not contain any sequence whose closure in $M^X$ is homeomorphic to $\beta\omega$.
\item $\overline{G}^{M^X}$ is a Rosenthal compactum.
\item $\vert \overline{G}^{M^X}\vert\leq \mathfrak{c}$.
\end{enumerate}
\ecor
\bpf
It suffices to apply Corollary \ref{resultado_1} and Lemma \ref{lema_K}.
%
%
\epf

\bcor\label{cor_1m}
Let $X$ be a compact space, $(M,d)$ be a metric space and $G$ a countable subset of $C(X,M)$ such that $\overline{G}^{M^X}$ is compact. If $\vert\overline{G}^{M^X}\vert\geq 2^{\mathfrak{c}}$, then there is a countable subset $L$ of $G$ such that its closure is canonically homeomorphic to  $\beta\omega$.
\ecor
\bpf
Use Corollary \ref{resultado_2}.
\epf

In the case where the metric space $M$ is $\C$ we have the following Corollary.


%

\bcor\label{resultado_2_C}
Let $X$ be a compact space and $G$ a uniformly and countable subset of $C(X,\C)$.
The following assertions are equivalent:
\begin{enumerate}[(a)]
\item $tg(\overline{G}^{\C^X})\leq \omega$.
\item $G$ does not contain any sequence whose closure in $\mathbb C^X$ is homeomorphic to $\beta\omega$.
\item $\overline{G}^{\C^X}$ is a Rosenthal compactum.
\item $\vert \overline{G}^{\C^X}\vert\leq \mathfrak{c}$.
\item $G$ does not contain a subsequence equivalent to the $\ell_1$ basis.
\end{enumerate}
\ecor
\bpf
Apply the complex version of Theorem \ref{talagrand} and Corollary \ref{resultado_2}.
\epf

\section{Dichotomy-type result for locally compact groups}
\mkp

Recall from the Introduction that for a locally compact group $G$ and if $\mathcal{H}_n\defi \C^n$ for $n=1,2, \ldots$;
$\mathcal H_{0}\defi l^2(\Z)$, then $Irr^C_n(G)$  denotes the set of irreducible unitary representations of $G$
on $\mathcal{H}_n$ (where it is assumed that every set $Irr^C_n(G)$ is equipped with the compact open topology),
and $Irr^C(G)=\bigsqcup\limits_{n\geq 0} Irr^C_n(G)$ (the disjoint topological sum).

The symbols $G^w$ (resp. $G^{w_C}$) designate the group $G$ equipped with the weak (group) topology generated by $Irr(G)$ (resp. $Irr^C(G)$).
As it was mentioned in the Introduction, if $G$ is abelian,
then the weak topology of $G$ coincides with the so-called Bohr topology associated to $G$.

\bdfn\label{I(G)}
We denote by $P(G)$ the set of continuous positive definite functions on $(G,\tau)$.
If $\sigma\in Irr(G)$ and $v\in \mathcal{H}^{\sigma}$, 
then the positive definite function:
$$\varphi:g\mapsto \langle\sigma (g)(v),v\rangle\text{,   }  g\in G$$
is called \textit{pure}, and the family of all such functions is denoted by $I(G)$. We also can define $I^C(G)$ the subset of $I(G)$ consisting of the elements whose irreducible representation is in $Irr^C(G)$. When $G$ is abelian, the set $I(G)$ coincides with the dual group $\widehat G$ of the group $G$.
\edfn
\mkp

The proof of the lemma below is straightforward.

\blem\label{top_debil_I}
Let $G$ be a locally compact group. Then:
\begin{enumerate}[(a)]
\item $G^{{\rm w}} = (G, {\rm w} (G,I(G)))$.
\item $G^{\w_C} = (G, \w(G,I^C(G)))$.
\end{enumerate}
\elem

\brem\label{Rem_top_debil_I}
We recall that $G^{{\rm w}} = G^{{\rm w_C}}$ if $G$ is a separable, metrizable, locally compact group.
\erem
\mkp

\bdfn\label{compactification wG}
Let $G$ be a locally compact group and consider the two following natural embeddings:
\begin{equation*}
\w :G\hookrightarrow \prod\limits_{\varphi\in I(G)}\overline{\varphi(G)}\hspace{1.5cm}  \hbox{and}\hspace{1.5cm}
\w_{C}:G\hookrightarrow \prod\limits_{\varphi\in I^C(G)}\varphi(G)
\end{equation*}
\begin{equation*}
\w (g)=(\varphi(g))_{\varphi\in I(G)} \hspace{3.8cm}\w_{C}(g)=(\varphi(g))_{\varphi\in I^C(G)}
\end{equation*}

\mkp

We define the \emph{weak compactification} $\w G$ (resp. \emph{$C$-weak compactification} $\w_CG$) of $G$ as the pair $(\w G,\w)$ (resp. $(\w_C G,\w_C)$),
where $\w G\defi \overline{\w (G)}$ (resp. $\w_CG\defi \overline{\w _C(G)}$).

This compactification has been previously considered in \cite{Cheng2011,Cheng2013} using different techniques. Also Akemann and Walter \cite{Akemann1972}
extended Pontryagin duality to non-abelian locally compact groups using the family of pure positive definite functions.
Again, in case $G$ is abelian, both compactifications, $(\w G,\w)$ and  $(\w_C G,\w)$, coincide with the
Bohr compactification of  $G$.

A better known compactification of a locally group $G$ which is closely related to $\w G$ is
defined as follows (cf. \cite{Eymard1964,Spronk2013}): let $\overline{\mathrm{B(G)}}^{\|\cdot\|_\infty}$ denote
the commutative $C^*$-algebra consisting of the uniform closure of the Fourier-Stieltjes algebra
of $G$. Here, the Fourier-Stieltjes algebra is defined as the matrix coefficients of
the unitary representations of $G$. Following \cite{Mayer1997} we call the spectrum $eG$ of
$\overline{\mathrm{B(G)}}^{\|\cdot \|_\infty}$ the
\emph{Eberlein compactification} of $G$.
Since the Eberlein compactification $eG$ is defined using the family of all continuous positive definite functions, it follows
that $\w G$ is a factor of $eG$ and, as a consequence,  inherits most of its properties. In particular,
$\w G$ is a compact involutive semitopological semigroup. 
\edfn
\mkp

The following definition was introduced by Hartman and Ryll-Nardzewski for abelian locally compact groups \cite{Hartman1964}.
Here, we extend it to arbitrary not necessarily abelian locally compact groups.

\bdfn
A subset $A$ of a locally compact group $G$ is an \emph{$I_0$ set} if every bounded complex (or real) valued function on $A$
can be extended to a continuous function on $\w G$. This definition extends the classic one, since when
$G$ is an abelian group, we have that ${\w G}=bG$, the so called \emph{Bohr compactification} of $G$
and $C(bG)\vert_G=AP(G)$ is the set of almost periodic functions on $G$.
\edfn

\brem\label{beta_I0}
Observe that if $(G,\tau)$ is a locally compact group and $A$ be a countably infinite subset of $G$, then $A$ is an $I_0$ set
if and only if $\overline A^{wG}$ is canonically homeomorphic to $\beta\omega$.
\erem
\mkp

The following Lemma  can be found in \cite[Section 14, Th.3]{Todorcevic1997}.

\blem\label{beta_Todorcevic}
Let $X$ be a compact space and $f:X\rightarrow \beta\omega$ a continuous and onto map. If $f^{-1}(n)$ is a singleton for all $n<\omega$ and $f^{-1}(\omega)$ is dense in $X$. Then $f$ is a homeomorphism.
\elem

\blem\label{w_ws}
Let $(G,\tau)$ be a {separable metric locally compact} group and $\lbrace g_n\rbrace_{n<\omega}$ be a sequence on $G$ such that $\overline{\lbrace g_n\rbrace_{n<\omega}}^{w_CG}\cong\beta\omega$, then $\overline{\lbrace g_n\rbrace_{n<\omega}}^{wG}\cong\beta\omega$.
\elem
\bpf
Let $\varphi:G^w\rightarrow G^{w_C}$ be the identity map, which is clearly a continuous group homomorphism and set
$\overline{\varphi}:wG\rightarrow w_CG$ the continuous extension of $\varphi$.
The result follows from Lemma \ref{beta_Todorcevic}.
\epf
\mkp

We now recall some known results about unitary representations of locally compact groups that are needed in the proof of our main result in this section.
One main point is the decomposition of unitary representations by direct integrals of irreducible unitary representations.
This was established by Mautner \cite{Mautner1950} following the ideas introduced by von Neuman in \cite{vonNeumann1949}.

\bthm[F. I. Mautner, \cite{Mautner1950}]\label{Th_Mautner}
For any representation $(\sigma,\mathcal H_\sigma)$ of a separable locally compact group $G$, there is a measure space
$(R,\mathcal R,r)$, a family $\{\sigma(p)\}$ of irreducible representations of $G$, which are associated
to each $p\in R$, and an isometry $U$ of $\mathcal H_\sigma$ such that
$$U\sigma U^{-1}=\int_R \sigma(p)d_r(p).$$
\ethm

\brem\label{Re_integral decomposition}
The proof of the above theorem given by Mautner assumes that the representation space $\mathcal H_\sigma$ is separable
but, subsequently, Segal \cite{Segal1951} removed this constraint. Furthermore, it is easily seen that
we can assume that $\sigma(p)$ belongs to $Irr^C(G)$ locally almost everywhere in the theorem above (cf. \cite{Ikeshoji1979a}).
\erem

A remarkable consequence of Theorem \ref{Th_Mautner}
is the following corollary about positive definite functions.

\bcor\label{Co_Mautner}
Every Haar-measurable positive definite function $\varphi$ on a separable locally compact group $G$ can be
expressed for all $g\in G$ outside a certain set of Haar-measure zero in the form
$$\varphi(g)= \int_R \varphi_p(g)d_r(p),$$

\noindent where $\varphi_p$ is a pure positive definite functions on $G$ for all $p\in R$.
\ecor
\mkp

The following proposition is contained in the proof of Lemma 3.2 of Bichteler \cite[pp. 586-587]{Bichteler1969}

\bprp\label{Pro_referee}
Let $G$ be a locally compact group. If $H$ is an open subgroup of $G$, then
each continuous irreducible representation of $H$
is the restriction of a continuous irreducible representation of $G$.
\eprp

\bdfn
Let $U$ be an open neighbourhood of the identity of a topological group $G$. We say that a sequence $\lbrace g_n\rbrace_{n<\omega}$ is \textit{$U$-discrete} if $g_nU\cap g_mU=\emptyset$ for all $n\neq m\in\omega$.
\edfn

\begin{proof}[\textbf{Proof of Theorem \ref{Teo_A}}]
Since $G$ is metric, we may assume  without loss of generality that the sequence is not $\lbrace g_n\rbrace_{n<\omega}$ is not $\tau$-precompact.
Otherwise, it would contain a $\tau$-convergent subsequence that, as a consequence, would be weakly convergent and
\emph{a fortiori} weakly Cauchy.

Thus, $\lbrace g_n\rbrace_{n<\omega}$ must contain  a subsequence that is $U_0$-discrete for some symmetric, relatively compact and open
neighbourhood of the identity $U_0$ in $G$. For simplicity's sake, we assume  without loss of generality that the whole sequence $\lbrace g_n\rbrace_{n<\omega}$ is
$U_0$-discrete.

Take the $\sigma$-compact, open subgroup $H\defi\langle \overline U_0\cup\lbrace g_n\rbrace_{n<\omega} \rangle $ of $G$.
Since $H$ is metric, $\sigma$-compact, it follows that $H$ is a Polish locally compact group. Consequently,
by \cite[Section 18.1.10]{Dixmier1964}, we have that  $Irr^C_m(H)$, equipped with the compact open topology,
is a Polish space for all $m\in \lbrace 0,1,2,\ldots\rbrace$.

The space $\mathcal{H}_m$ being separable, for each $m\in\lbrace 0,1,2,\ldots\rbrace$, there exists a countable subset
$D_m\defi\lbrace v^m_n\rbrace_{n<\omega}$ that is dense in the unit ball of $\mathcal{H}_m$ (therefore,
the linear subspace generated by $D_m$ will be dense in $\mathcal{H}_m)$.
Fix $m\in\lbrace 0,1,2,\ldots\rbrace$ and 
let $\bD$ denote the closed unit disk in $\C$.
We have that $<\sigma(g)(v^m_n),v^m_n>\in \bD$ for all $\sigma\in Irr_m^C(H)$, $g\in H$ and $n<\omega$.

For each $m\in\lbrace 0,1,2,\ldots\rbrace$, let $\alpha_m:H\rightarrow C_p(Irr^C_m(H),\bD^\omega)$ be the continuous and injective map defined by $\alpha_m(h)(\sigma)\defi (<\sigma(h)(v^m_n),v^m_n>)_{n\in\omega}$ for all $h\in H$ and $\sigma\in Irr_m^C(H)$.
Since $\bD^\omega$ is a compact metric space, it follows that
$$\overline{\lbrace \alpha_m(g_n)\rbrace_{n<\omega}}^{\,(\bD^\omega )^{Irr^C_m(H)}}$$ is compact for all $m\in\lbrace 0,1,2,\ldots\rbrace$.

Now, we successively apply Corollary \ref{cor_Rosenthal_metrico} for each $m\in\lbrace 0,1,2,\ldots \rbrace$ as follows.

For $m=0$, $\lbrace \alpha_0(g_n)\rbrace_{n<\omega}$ contains either a pointwise Cauchy subsequence
or a subsequence whose closure in $(\bD^\omega)^{Irr_0^C(H)}$ is canonically homeomorphic to $\beta\omega$.

If there is a Cauchy subsequence $\lbrace \alpha_0(g_{n^0_i})\rbrace_{i<\omega}$, then we go on to the case $m=1$.
That is  $\lbrace \alpha_1(g_{n^0_i})\rbrace_{i<\omega}$ contains either a pointwise Cauchy subsequence or a subsequence whose closure in
$(\bD^\omega)^{Irr_1^C(H)}$ is canonically homeomorphic to $\beta\omega$. If there is a Cauchy subsequence
$\lbrace \alpha_1(g_{n^1_i})\rbrace_{i<\omega}$ we go on to the case $m=2$, and so forth.

Assume that we can find a pointwise Cauchy subsequence in each step and take the diagonal subsequence $\lbrace g_{n^i_i}\rbrace_{i<\omega }$.
We have that $\lbrace \alpha_m(g_{n^i_i})\rbrace_{i<\omega }$ is pointwise Cauchy for each $m\in\lbrace 0,1,2,\ldots\rbrace$.
We claim that the the subsequence $\lbrace g_{n^i_i}\rbrace_{i<\omega }$ is Cauchy in the weak topology of $G$.


Indeed, take an arbitrary element $\varphi\in I^C(H)$, then there is $t\in \lbrace 0,1,2,\ldots\rbrace$, $\sigma\in Irr^C_t(H)$ and $v\in \mathcal{H}_t$ such that
$\varphi(h)=<\sigma(h)(v),v>$ for all $h\in H$, where we may assume that $\|v\|\leq 1$  without loss of generality.

Let $\epsilon >0$ be an arbitrary positive real number. By the density of $D_t$, there is $u\in D_t$ such that
$\|u-v\|<\epsilon/6$.

For every $h\in H$, we have

\begin{align*}
\vert <\sigma(h)(v),v> - <\sigma(h)(u),u> \vert &=
\vert <\sigma(h)(v),v> - <\sigma(h)(u),v> +\\
&\quad \ \ <\sigma(h)(u),v> - <\sigma(h)(u),u> \vert \\ &\leq
\vert <\sigma(h)(v-u),v>\vert + \vert<\sigma(h)(u),v-u> \vert\\ &\leq
2\|v-u\|<\epsilon/3.
\end{align*}

On the other hand, we have that $\lbrace \alpha_t(g_{n^i_i})\rbrace_{i<\omega }$ is a pointwise
Cauchy sequence in $(\bD^\omega )^{Irr^C_t(H)}$. Thus,
from the definition of $\alpha_t$ and, since $u\in D_t$,
it follows that $$\lbrace<\sigma(g_{n^i_i})(u),u> \rbrace _{i<\omega}$$ is a Cauchy sequence in $\bD$.
Hence, there is $i_0<\omega$ such that

$$\vert <\sigma(g_{n^i_i})(u),u> - <\sigma(g_{n^j_j})(u),u> \vert < \epsilon/3\ \text{for all}\ i,j\geq i_0.$$
\mkp

This yields
\begin{align*}
\vert <\sigma(g_{n^i_i})(v),v> - <\sigma(g_{n^j_j})(v),v> \vert &\leq
\vert <\sigma(g_{n^i_i})(v),v> - <\sigma(g_{n^i_i})(u),u> \vert\\ &+
\vert <\sigma(g_{n^i_i})(u),u> - <\sigma(g_{n^j_j})(u),u> \vert\\ &+
\vert <\sigma(g_{n^j_j})(u),u> - <\sigma(g_{n^j_j})(v),v> \vert\\ &<
3\, \epsilon/3=\epsilon.
\end{align*}
\mkp

We conclude that $\lbrace <\sigma(g_{n^i_i})(v),v> \rbrace _{i<\omega}=\lbrace \varphi(g_{n^i_i})\rbrace _{i<\omega}$
is a Cauchy sequence in $\bD$ for all
$\varphi\in I^C(H)$. Since $H$ is a locally compact Polish group,
we have that $G^{{\rm w}} = G^{{\rm w_C}}$ by Remark \ref{Rem_top_debil_I}. As a consequence, it follows that
which proves that $\lbrace g_{n^i_i}\rbrace_{i<\omega }$ is weakly Cauchy in $H$.
We must now verify that $\lbrace g_{n^i_i}\rbrace_{i<\omega }$ is weakly Cauchy in $G$.

In order to do so, take a map $\psi\in I(G)$. Since $H$ is separable, by Corollary \ref{Co_Mautner},
there is a measure space $(R,\mathcal R,r)$, a family  $\{\psi_p\}$ of pure positive definite functions on $H$, which are associated
to each $p\in R$, such that
$$\psi(h)= \int_R \psi_p(h)d_r(p)\ \hbox{for all}\ h\in H.$$

Therefore

$$\psi(g_{n^i_i})= \int_R \psi_p(g_{n^i_i})d_r(p))\ \hbox{for all}\ i<\omega.$$

Now, for each $i<\omega$, consider the map $f_i$ on $R$ by $f_i(p)\defi \psi_p(g_{n^i_i})$.
Then $f_i$ is integrable on $R$ and, since  $\lbrace g_{n^i_i}\rbrace_{i<\omega }$ is weakly Cauchy in $H$,
it follows that $\{f_i\}$ is a pointwise Cauchy sequence on $R$. Furthermore, if
$\psi_p(h)=<\sigma_p(h)[v_p],v_p>$ for some $\sigma_p\in Irr(H)$ and $v_p\in\mathcal H_{\sigma_p}$,
it follows that $$|f_i(p)|=|\psi_p(g_{n^i_i})|=|<\sigma_p(g_{n^i_i})[v_p],v_p>|\leq \|v_p\|^2.$$
Thus  defining $f$ on $R$ as the pointwise limit of $\{f_i\}$,
we are in position to apply Lebesgue's dominated convergence theorem in order to obtain that
$$\int_R f(p)d_r(p)=\lim\limits_{i\rightarrow \infty} \int_R \psi_p(g_{n^i_i})d_r(p)=\lim\limits_{i\rightarrow \infty} \psi(g_{n^i_i}).$$

In other words, the sequence $\{\psi(g_{n^i_i})\}$ converges and, therefore, is Cauchy for all $\psi\in I(G)$.
Hence $\{g_{n^i_i}\}$ is weakly Cauchy in $G$ and we are done.

Suppose now that there exists an index $m_0\in \lbrace 0,1,2,\ldots,\infty\rbrace$ such that
$\lbrace \alpha_{m_0}(g_{n^{m_0}_{\,i}})\rbrace_{i<\omega}$ contains a subsequence $\lbrace \alpha_{m_0}(g_{n(j)})\rbrace_{j<\omega}$
whose closure in $(\bD^\omega)^{Irr_{m_0}^C(H)}$ is homeomorphic to $\beta\omega$.
Applying Lemma \ref{beta_Todorcevic}, we know that $\overline{\lbrace g_{n(j)}\rbrace_{j<\omega}}^{w_CH}\cong \beta \omega$.
Consequently, by Lemma \ref{w_ws}, we obtain that $\overline{\lbrace g_{n(j)}\rbrace_{j<\omega}}^{wH}\cong\beta\omega$.

On the other hand, by Proposition \ref{Pro_referee}, 
we have that the irreducible representations of $H$
are the restrictions of irreducible representations of $G$, which implies that the identity map
$id:(H,w(G,I(G))_{|H})\rightarrow (H,w(H,I(H)))$
is a continuous group isomorphism that can be extended canonically to a homeomorphism between their associated compactifications
$\overline{id}:\overline{H}^{wG}\rightarrow wH$.
By Lemma \ref{beta_Todorcevic} again, we obtain that
 $\overline{\lbrace g_{n(j)}\rbrace_{j<\omega}}^{wG}\cong \beta \omega$.
 Thus $\lbrace g_{n(j)}\rbrace_{j<\omega}$ is an $I_0$ set, which completes the proof.
\end{proof}

\brem
Theorem \ref{Teo_A} fails if we try to extend it to every locally compact group or even to every compact group. Indeed,
Fedor\v{c}uk \cite{Fedorcuk1977} has proved that the existence of a compact space $K$ of cardinality $\mathfrak{c}$ without convergent sequences is compatible with \textbf{ZFC}.
If we take the Bohr compactification of the free abelian group generated by $K$, then every sequence contained in $K$ does not fulfil any of the two choices
established in Rosenthal's dichotomy.
\erem

\section{Interpolation sets}

Hartman and Ryll-Nardzewski \cite{Hartman1964} proved that every abelian locally compact group contains an $I_0$ set.
This result was improved in \cite{Galindo_Hernandez1999}, where it was proved that
every non-precompact subset of an abelian locally compact group contains an $I_0$ set.
These sort of results do not hold for general locally compact groups unfortunately. Indeed,
the Eberlein compactification of the group $SL_2(R)$ coincides with its one-point compactification,
which means that each continuous positive definite function on $SL_2(R)$
converges at infinity (see \cite{Chou1980}). Therefore, for this group, only the first case of the dichotomy
result in Theorem \ref{Teo_A} holds. If we search for interpolation sets, some extra conditions have to be assumed.

In this section we explore the application of the results in the previous sections
in the study of interpolation sets in locally compact groups. First, we need the following
result that was established by Ernest \cite{Ernest1971} (cf. \cite{Ikeshoji1979})
for separable metric locally compact groups and convergent sequences and
Subsequently extended for locally compact groups and compact subsets by Hughes \cite{Hughes1972} .

\bprp\emph{(J. Ernest, J.R. Hughes)}\label{compact_3}
Let $(G,\tau)$ be a locally compact group. Then $(G,\tau)$ and $G^w$ contain the same compact subsets.
\eprp
\mkp

In some special cases, Hughes' result implies the convergence of weakly Cauchy sequences.

\bprp\label{compact_4}
Let $(G,\tau)$ be a locally compact group and suppose that $\{g_{n}\}_{n<\omega}$ is a Cauchy sequence in $G^w$.
If $\overline{\{g_{n}\}_{n<\omega}}^{\, \w G}\subseteq inv(\w G)$, then $\{g_{n}\}_{n<\omega}$ is $\tau$-convergent in $G$.
\eprp
\bpf
Assume that $\{g_{n}\}_{n<\omega}$ is a Cauchy sequence in $G^w$. First, we verify that the sequence is
a precompact subset of $(G,\tau)$.

Indeed, we have that $\{g_{n}\}_{n<\omega}$ converges to some element $p\in inv(\w G)$.
If $\{g_{n}\}_{n<\omega}$ were not precompact in $(G,\tau)$, there would be a neighbourhood of the
neutral element $U$ and a subsequence $\{g_{n(m)}\}_{m<\omega}$ such that
$g_{n(m)}^{-1}\cdot g_{n(l)}\notin U$ for each $m,l<\omega$ with $m\neq l$.
On the other hand, the sequence $\{g_{n(m)}^{-1}\cdot g_{n(m+1)}\}_{m<\omega}$
converges to $p^{-1}p$, the neutral element in $G^w$.
This takes us to a contradiction because, by Proposition \ref{compact_3},
it follows that $\{g_{n(m)}^{-1}\cdot g_{n(m+1)}\}_{m<\omega}$ must also converge to the neutral element in  $(G,\tau)$.

Therefore, the sequence $\{g_{n}\}_{n<\omega}$ is a precompact subset of $(G,\tau)$. This
implies that $p\in G$ and we are done.
\epf
\mkp

\blem\label{Teo_C}
Let $(G,\tau)$ be a locally compact group and let $B$ be a non-precompact subset of $G$ such that
$\overline B^{\, \w G}\subseteq inv(\w G)$.
Then there exist a open subgroup $H$ of $G$, a compact and normal subgroup $K$ of $H$, a quotient map
$p:H\rightarrow H/K$ and a sequence $\lbrace g_n\rbrace_{n<\omega}\subseteq B\cap H$ such that $H/K$ is a
Polish group and $\overline{p(\lbrace g_n\rbrace_{n<\omega})}^{w H/K}\cong \beta\omega$.
\elem
\bpf
Since $B$ is non-precompact there exists an open, symmetric and relatively compact neighbourhood of the identity $U$ in $G$
such that $B$ contains a $U$-discrete sequence $\lbrace g_n\rbrace_{n<\omega}$.

Consider the subgroup $H\defi < \overline U\cup\lbrace g_n\rbrace_{n<\omega} >$, which $\sigma$-compact and open in $G$.
By Kakutani-Kodaira's theorem, there exists a normal, compact $K$ of $H$ such that $K\subseteq U$ and $H/K$ is metrizable, and consequently Polish.
Let $p:H\rightarrow H/K$ be the quotient map and let $\overline p :\w H\rightarrow \w H/K$ denote the canonical extension
to the weak compactifications. Therefore, we have that $\overline p(inv(\w H)\subseteq inv(\w H/K)$.
Furthermore, since $\overline H^{\, \w G}$ is canonically homeomorphic to $\w H$, it follows that
$\overline{\lbrace g_n\rbrace_{n<\omega}}^{\, \w H}\subseteq inv(\w H)$. Hence
$\overline{\lbrace p(g_n)\rbrace_{n<\omega}}^{\, \w H/K}\subseteq inv(\w H/K)$. Thus, we are in position of applying Proposition \ref{compact_4}.

Assume that there is a weakly Cauchy subsequence $\lbrace p(g_s)\rbrace_{s<\omega}$ in $(H/K)^\w$,
which would be $\tau/K$-convergent by Proposition \ref{compact_4}.
Then by a theorem of Varopoulos \cite{Varopoulos1964}, the sequence $\lbrace p(g_s)\rbrace_{s<\omega}$ could be lifted to
a sequence $\lbrace x_s\rbrace_{s<\omega}\subseteq H$ converging to some point $x_0\in H$. This would entail that $ x_s^{-1}g_s\in K$ for all $s\in\omega$.
Thus the sequence $\lbrace g_s\rbrace_{s<\omega}$ would be contained in the compact subset
$(\lbrace x_s\rbrace_{s<\omega}\cup\lbrace x_0\rbrace) K$, which is a contradiction since $\lbrace g_n\rbrace_{n<\omega}$ was supposed
to be $U$-discrete.
This contradiction completes the proof.
\epf
\begin{proof}[\textbf{Proof of Theorem \ref{Cor_C}}]
Applying Lemmata \ref{beta_Todorcevic} and \ref{Teo_C}, we obtain an $I_0$ set.
On the other hand, according to \cite[Proposition II.4.6]{Ruppert1984},
any intrinsic group at an idempotent in a semi-topological semi-group is complete with respect to the two-sided uniformity,
which means that $inv(\w G)$ is complete in the two-sided weak uniformity. As a consequence, our $I_0$ set is precompact
in the weak (group) completion of $G$.
\end{proof}
\mkp

\bcor\label{Cor_vD}
Let $G$ be a discrete group and let $\lbrace g_n\rbrace_{n<\omega}$ be an infinite sequence in $G$.
If $\overline{\lbrace (g_n)\rbrace_{n<\omega}}^{\, \w G}\subseteq inv(\w G)$, then $\lbrace g_n\rbrace_{n<\omega}$
contains an infinite $I_0$ set.
\ecor
\mkp

\brem
In case the group $G$ is abelian, Corollary \ref{Cor_vD} is a variant of van Douwen's Theorem \ref{douw}.
\erem

We now look at \emph{Sidon sets}, a well known family of interpolation sets in harmonic analysis.
In fact, there are several definitions for the notion of Sidon sets for nonabelian groups, although
all coincide for amenable groups (see \cite{Figa-Talamanca1977}). Here, we are interested in the following:
We say that a subset $E$ of a locally compact group $G$ is a \emph{(weak) Sidon set} if every bounded function
can be interpolated by a continuous function defined on the Eberlein compactification $eG$.
This is a weaker property than the classical notion of \textit{Sidon set}
in general (see \cite{Picardello1973}) but both notions coincide for amenable groups.

\begin{proof}[\textbf{Proof of Theorem \ref{Teo_Sidon}}]
Since every $I_0$ set is automatically weak Sidon, it suffices to apply Theorem \ref{Cor_C}.
\end{proof}
\mkp

We notice that the following question still remains open (see \cite{Lopez1975} and \cite[p. 57]{Figa-Talamanca1977}).

\bqtn[Fig\`{a}-Talamanca, 1977]\label{FT}
Do infinite weak Sidon exist in every infinite discrete non-abelian group?
\eqtn
\mkp

Theorem \ref{Teo_Sidon} and Corollary \ref{Cor_vD} provide
sufficient conditions for the existence of weak Sidon sets. 


\bdfn\label{independent}
Let $G$ be a group. A sequence $\lbrace x_n\rbrace_{n<\omega}\subseteq G$ is called
\textit{independent} if for every $n_0<\omega$ the element $x_{n_0}\notin \langle\lbrace x_n\rbrace_{n\in \omega\setminus n_0}\rangle$.
A group $G$ is called \emph{locally finite}
if every finite subset of the group generates a finite subgroup.
The group $G$ is \emph{residually finite} if for every
non-identity element $g$ of $G$ there exists a normal subgroup $N$
of finite index in $G$ such that $g\notin N$. Finally, the group
$G$ is called an $FC$-group if every conjugacy class of $G$ is
finite.
\edfn

\bprp\label{prop_indep_sidon}
Every independent sequence in a discrete group $G$ is a weak Sidon set.
\eprp
\bpf
Let $E=\lbrace x_n\rbrace_{n<\omega}$ be an independent sequence in $G$. It will suffice to show that $\overline{E}^{eG}$ is homeomorphic to $\beta\omega$ or,
equivalently, that every pair of disjoint subsets in $E$ have disjoint closures in $eG$.

Indeed, for any pair $A,B$ of arbitrary disjoint subsets of $\omega$, set $X_A\defi \lbrace x_n\rbrace_{n\in A}$,
and $X_B\defi \lbrace x_n\rbrace_{n\in B}\subseteq G\setminus \langle X_A\rangle$.
Since $\langle X_A\rangle$ is an open subgroup of $G$, the positive definite function $h$ defined by 
${h}(x)= 1$ if $x\in \langle X_A\rangle$ and ${h}(x)=0$ if $x\in G\setminus \langle X_A\rangle$ is continuous
(see \cite[32.43(a)]{Hewitt1970}). Thus $\overline{X_A}^{\, eG}\cap \overline{X_B}^{\, eG}=\emptyset$, which completes the proof.
\epf

\bcor\label{cor_indep_units}
Let $G$ be a discrete group and let $\lbrace x_n\rbrace_{n<\omega}$ be a sequence in $G$. If the sequence contains either an independent subsequence
or its $\w G$-closure is contained in $inv(\w G)$, then the sequence contains a weak Sidon set.
\ecor

\brem
Note that if $G$ is a discrete $FC$-group we can find an independent sequence within every infinite subset of $G$.
Therefore, every infinite subset of a discrete $FC$-group contains a weak Sidon set \cite{Michele1976}.
\erem

We finish this section with an example of a sequence that is a weak Sidon set and converges to the neutral element in the Bohr topology.

\bprp\label{le_bohr5}
Let $G$ be a residually finite, locally finite group that is not abelian by finite.
Then $G$ contains a sequence that is weak Sidon and converges in the Bohr topology.
\eprp
\begin{proof}
Let $a_1$ be a non trivial element in $G'$. Then, there is a
finite non Abelian subgroup $B_1$ of $G$ with $a_1\in B'_1$.
Suppose we have defined $L=\sum_{i=1}^nB_i$. Since $L$ is finite
and $G$ is residually finite, we can get $m\in \N$ and $\sigma\in
Rep_m(G)$ such that $\sigma_{|L}$ is faithful and $\sigma(G)$ is finite. Set
$N= ker \sigma$. If $N$ were Abelian, it would follow that $G$ is
Abelian by finite, which is impossible. Thus, we may assume  without loss of generality
that $N'\neq \{1 \}$ and we can replace $G$ by $N$ in order to
obtain a finite non Abelian subgroup $B_{n+1}$ of $N$ and a non
trivial element $a_{n+1}\in B'_{n+1}$. Using an inductive
argument, we obtain $B=\sum_{i=1}^{\infty} B_i$, a subgroup of
$G$, such that $B'_i\neq \{1 \}$, for all $i\in \N$. Under such
circumstances, it is known that the sequence $\{a_n\}$ is convergent
in the Bohr topology (cf. \cite[Cor. 3.10]{Hernandez2006}). On the other hand,
the sequence $\{a_n\}$ is an independent set by definition and, therefore,
a weak Sidon subset in $G$. This completes the proof.
\end{proof}

Remark that a positive answer to the following question will also solve Question \ref{FT}.

\bprb
Does every discrete group contain an infinite $I_0$ set?
\eprb

\section{Property of strongly respecting compactness}

One of the first papers dealing with the weak topology of locally compact groups was given by Glicksberg in \cite{Glicksberg1962}, where he proved that every weakly compact subset of a locally compact abelian group must be compact. Using the terminology introduced by Trigos-Arrieta in \cite{Trigos-Arrieta1991}, it is said that locally compact abelian groups
\emph{respect compactness}.
Along this line, Comfort, Trigos-Arrieta and Wu \cite{Comfort1993} established the following substantial extension of Glicksberg's theorem:
Let $G$ be a locally compact abelian group with Bohr compactification $bG$ and let $N$ be a metrizable closed subgroup of $bG$.
Then a subset $A$ of $G$ satisfies that $AN\cap G$ is compact in $G$ if and only if $AN$ is compact in $bG$.
Motivated by this result, the authors defined a group that \textit{strongly respects compactness} to be a group that satisfies the assertion above.
They also raised the question of characterizing such groups. This question was consider in \cite{Galindo_Hernandez2004} but looking at the finite dimensional, irreducible representations of a locally compact group. Here we consider it in complete generality.
First, we need the following definition that extends the notion of strongly respecting compactness
to groups that are not necessarily abelian equipped with the weak topology.

\bdfn
We say that a locally compact group $G$ \textit{strongly respects compactness} if for any closed metrizable subgroup $N$ of $inv(\w G)$,
a subset $A$ of $G$ satisfies that $AN\cap G$ is compact in $G$ if and only if $AN$ is compact in $\w G$.
\edfn

The main result of this section states that every locally compact group strongly respects compactness. 

Let $G$ be a locally compact group and let $H$ be an open subgroup of $G$. We have discussed in the proof of Theorem \ref{Teo_A}
that the identity map
$id:(H,\w(G,I(G))_{|H})\rightarrow H^\w $ is a continuous isomorphism that can be extended to a continuous map
$\overline{id}:\overline{H}^{\w G}\rightarrow \w H$ defined on their respective compactifications.
On the other hand, using that $\w G$ is factor of the Eberlein compactification $eG$, it
follows that $\w G$ is a compact involutive semitopological semigroup for every locally compact group $G$.
Taking this fact into account, the following lemma is easily verified.

\blem\label{Le_Semigroup}
Let $G$ and $H$ be locally compact groups and let $h\colon G\to H$  be a continuous homomorphism.
Then there is a canonical continuous extension $\overline{h}:\w G\rightarrow \w H$ such that
for every $p,q$ in $\w G$, we have $\overline{h}(pq)=\overline{h}(p)\overline{h}(q)$.
\elem

The next result is a variation of \cite[Lemma 3.6]{Galindo_Hernandez2004}.

\blem\label{Lema_1}
Let $G$ be a locally compact group, $H$ an open subgroup of $G$, $A$ a subset of $G$, and let $N$ be a subgroup of $inv(\w G)$, containing the identity,
such that $AN$ is compact in $\w G$. If $F$ is an arbitrary subset of $AN\cap H$, then there exists $A_0\subseteq A$ with $\vert A_0\vert\leq \vert N\vert$ such that $$\overline{F}^{\w H}\subseteq \overline{(FF^{-1})}^{H^\w }\cdot\overline{id}(A_0 N).$$
\elem
\bpf
We first verify that $\overline{F}^{\w H}\subseteq \overline{id}( AN\cap \overline{H}^{\w G})$.

Indeed,
since $AN\cap H\subseteq AN\cap \overline{H}^{\w G}$ and $AN\cap\overline{H}^{\w G}$ is compact, we have $AN\cap H\subseteq \overline{id}( AN\cap \overline{H}^{\w G})$ and,
as a consequence, it follows that $\overline{AN\cap H}^{\w H}\subseteq \overline{id}( AN\cap \overline{H}^{\w G})$. Hence $\overline{F}^{\w H}\subseteq \overline{id}( AN\cap \overline{H}^{\w G})$ and $\overline{F}^{\w H}$ is compact.

For any $x\in N$ such that $\overline{F}^{\w H}\cap\overline{id}(Ax)\neq\emptyset$, pick $a_x\in A$ with $\overline{id}(a_xx)\in \overline{F}^{\w H}$.
We define $A_0\defi \lbrace a_x\in A :x\in N \ \hbox{and}\ \overline{id}(a_xx)\in \overline{F}^{\w H} \rbrace$. We have $A_0\subseteq A$ and $\vert A_0\vert\leq \vert N\vert$.

Pick an arbitrary point $b\in \overline{F}^{\w H}$. Since  $\overline{F}^{\w H}\subseteq \overline{id}( AN\cap \overline{H}^{\w G})$ we can find
$a\in A$ and $y\in N$ such that $b=\overline{id}(ay)$. Set $b'=\overline{id}(a_yy)\in \overline{F}^{\w H}$.
Then $bb'^{-1}=\overline{id}(ay)\overline{id}(a_yy)^{-1}\in \overline{F}^{\w H}\overline{F^{-1}}^{\w H}=\overline{FF^{-1}}^{\w H}$.
Observe also that, by Lemma \ref{Le_Semigroup}, we have $bb'^{-1}=\overline{id}(ayy^{-1}a_y^{-1})=\overline{id}(aa_y^{-1})=aa_y^{-1}\in \w H\cap G=H$.
Therefore $bb'^{-1}\in \overline{FF^{-1}}^{\w H}\cap H$.
Since $H$ is an open subgroup, by \cite[Cor. 14.2]{Hughes1972}, we deduce that $bb'^{-1}\in \overline{FF^{-1}}^{\w H}\cap H=\overline{FF^{-1}}^{H^\w }$.
Thus $b=bb'^{-1}b'\in \overline{FF^{-1}}^{H^\w }\cdot \overline{id}(A_0N)$\ and we are done.
\epf
\mkp

\begin{proof}[\textbf{Proof of Theorem \ref{Teo_SRC}}]
Let $N$ be a metrizable subgroup of $inv(\w G)$ and let $A$ a subset of $G$ such that $AN$ is compact in $\w G$. Since $AN\cap G$ is closed in $G$ it suffices to see that it is precompact. By reduction to absurd, assume that $AN\cap G$ is non-precompact. By Theorem \ref{Teo_C} there exist an open subgroup $H$ of $G$, a compact and normal subgroup $K$ of $H$, a quotient map $p:H\rightarrow H/K$ and a sequence $F\subseteq AN\cap H$ such that $H/K$ is a Polish group and $\overline{p(F)}^{\w  H/K}\cong \beta\omega$. Thus, $\vert \overline{p(F)}^{\w  H/K}\vert\geq 2^{\mathfrak{c}}$.

By Lemma \ref{Lema_1} there is $A_0\subseteq A$ with $\vert A_0\vert\leq \vert N\vert$ such that $\overline{F}^{\w H}\subseteq \overline{(FF^{-1})}^{H^\w }\overline{id}(A_0 N)$.
Since $\vert \overline{id}(A_0 N)\vert\leq \mathfrak{c}$ we can enumerate it as $\lbrace a_{\alpha}\rbrace_{\alpha<\mathfrak{c}}$. Therefore, we can write $\overline{F}^{\w H}\subseteq\bigcup\limits_{\alpha<\mathfrak{c}} \overline{FF^{-1}}^{H^\w }a_{\alpha}$.

Let $\overline p:\w H\rightarrow \w H/K$ be the canonical extension of $p$ to the respective compactifications of $H$ and $H/K$.
Using Lemma \ref{Le_Semigroup}, for each $z\in \w H$ consider the map $T_z$ defined on $\w  H/K$ by
$T_z(\overline p(x))=\overline p(xz)=\overline p(x)\overline p(z)$ for all $x\in \w H$. Hence, from the previous inclusion we obtain that $\overline{p(F)}^{\w  H/K}=\overline{p}(\overline{F}^{\w H})\subseteq\bigcup\limits_{\alpha<\mathfrak{c}} T_{a_{\alpha}}(\overline{p(FF^{-1})}^{H^\w /K})$.
Since the topology of $H^\w /K$ is finer than that of $(H/K)^\w $ we have that $\overline{p(FF^{-1})}^{H^\w /K}\subseteq \overline{p(FF^{-1})}^{(H/K)^\w }$.
Furthermore $\vert \overline{p(FF^{-1})}^{(H/K)^\w }\vert \leq \mathfrak{c}$ because $H/K$ is a Polish space.
Therefore, $\vert \overline{p(F)}^{\w  H/K}\vert\leq \mathfrak{c}$. This is a contradiction that completes the proof.
\end{proof}


\section{Acknowledgments}
We are very grateful to the referee for several valuable comments and suggestions that helped to improve the presentation of this paper.

\bibliography{BiblioResultados}
\bibliographystyle{plain}

non-precompact subset
\end{document}